\documentclass{amsart}
\usepackage{amssymb}
\newcommand{\cC}{{\mathcal C}}

\newcommand{\cM}{{\mathcal M}}
\newcommand{\be}{{\bf 1}}

\newcommand{\BZ}{{\mathbb Z}}

\newcommand{\sq}{$\square$}

\newcommand{\Hom}{\mbox{Hom}}

\newcommand{\Id}{\mbox{Id}}

\title{Fusion categories of rank 2}
\author{Viktor Ostrik}
\email{ostrik@math.mit.edu}
\address{Department of Mathematics, MIT, 77 Massachusetts Ave., Cambridge,
MA 02139}
\thanks{The author was partially supported by NSF grant DMS-0098830}
\date{March 2002}
\begin{document}
\begin{abstract} We classify semisimple rigid monoidal categories with two
isomorphism classes of simple objects over the field of complex numbers.
In the appendix written by P.~Etingof it is proved that the number of 
semisimple Hopf algebras with a given finite number of irreducible
representations is finite. 
\end{abstract} 
\maketitle
\section{Introduction}
Let $k$ be an algebraically closed field of characteristic 0.
A fusion category $\cC$ over $k$ is a $k-$linear semisimple rigid monoidal
category with finitely many simple objects and finite dimensional spaces
of morphisms, such that the endomorphism algebra of the neutral object is $k$,
see \cite{ENO}. 

The problem of classification of all fusion categories seems to be very 
difficult. A basic invariant of a fusion category $\cC$ is its Grothendieck 
ring $K(\cC)$ (physically, fusion rules algebra or fusion ring) which is a 
unital based ring in the sense of \cite{O}. It is
known that for a given based ring $K$ there are only finitely many fusion
categories $\cC$ with $K(\cC)=K$ (this statement is known as ``Ocneanu 
rigidity'', see \cite{W} and \cite{ENO}). There is a natural problem for a 
given ``interesting'' based ring $K$ to find all fusion categories $\cC$ with
$K(\cC)=K$ (such categories are called ``categorifications'' of $K$). 
This problem was considered first probably by G.~Moore and N.~Seiberg in
\cite{MS}, they considered the case when all objects of $\cC$ are invertible
and the case of Yang-Lee fusion rules (see below). In his thesis T.~Kerler
completely classified fusion categories with fusion rings isomorphic to
the fusion ring of the category of integrable $\widehat{sl}_2-$modules at
a positive integer level, see \cite{Ke}, and later D.~Kazhdan and H.~Wenzl 
generalized this to the case of $\widehat{sl}_n$, see \cite{KW}. In \cite{TY}
D.~Tambara and S.~Yamagami considered another big class of examples, the
so-called fusion rules of self-duality for finite abelian groups.  

The rank of $K(\cC)$ over $\BZ$ or equivalently the number of isomorphism 
classes of simple objects in $\cC$ is called the rank of the category $\cC$.
There is only one fusion category of rank 1 --- the category of vector spaces
over $k$. In this note we will study all fusion categories of rank 2. Let
$\be, X$ be the simple objects of such category $\cC$ (here $\be$ is the
unit object). It is clear that possible fusion rules for $\cC$
are completely determined by the number $n\in \BZ_{\ge 0}$  from the equation
$X\otimes X=\be \oplus nX$. Let $K_n$ denote the fusion ring corresponding to
the number $n$. The fusion ring $K_1$ is called the Yang-Lee fusion rules and
is well known in the conformal field theory.

{\bf Main Theorem.} There are just 4 fusion categories of rank 2. For two of
them $K(\cC)=K_0$ and for other two of them $K(\cC)=K_1$.

The classification of fusion categories $\cC$ with $K(\cC)=K_0$ or $K(\cC)=
K_1$ is due to Moore and Seiberg \cite{MS}. The only new part of this Theorem
is the fact that the fusion rings $K_n, n\ge 2$ admit no categorification. 

Our result suggests that the answer to the following question is of interest:

{\bf Question.} Is it true that there are only finitely many fusion 
categories of a given finite rank?

As a supporting evidence for positive answer to this question recall that the
number of finite groups with a given number of irreducible representations is
finite, see \cite{L}. Moreover, P.~Etingof proved that the number of finite
dimensional semisimple Hopf algebras with a given number of irreducible
representations is finite, see Appendix to this note. 

{\bf Remark.} Our Main Theorem is not true without rigidity assumption on the
category $\cC$. See \cite{V} for an example of a semisimple 
bialgebra with two representations and fusion rules $X\otimes X=2X$.

I am much obliged to Tania Chmutova, Pavel Etingof and Dmitri Nikshych for
useful discussions; many thanks are due to Pavel Etingof for writing the
Appendix to this note. I am grateful to Richard Stanley for 
providing reference \cite{L}.

\section{Proof of the Main Theorem}
It appears that direct methods developed in \cite{TY} are very difficult to
apply in a case when some fusion coefficients are greater than 1. In 
particular I could not study even the case of the fusion ring $K_2$ using
these methods. So we are using another approach. Here is an outline of our
proof. Let $\cC$ be a fusion category with $K(\cC)=K_n$.
First we study the Drinfeld double of $\cC$ and show that $\cC$ is braided. 
Then we show that the category $\cC$ is automatically ribbon. It is easy to
see then that the category $\cC$ is modular (unless $n=0$). Now the standard
identities from the theory of modular categories give us a contradiction.

\subsection{The category $\cC$ is braided} Let $Z(\cC)$ be the Drinfeld 
center of the category $\cC$, see e.g. \cite{M2}. It is known that the
category $Z(\cC)$ is semisimple, see \cite{ENO}. Let $\cC^{op}$ denotes 
the opposite category to $\cC$ and let $\cC\boxtimes \cC^{op}$ be the 
external product of the categories $\cC$ and $\cC^{op}$, see e.g. \cite{M2}.  

{\bf Lemma 2.1.} The category $Z(\cC)$ has 4 simple objects $\be, X_1, X_2,
Y$. Under the forgetful functor $F: Z(\cC)\to \cC$ one has $F(\be)=\be,
F(X_1)=F(X_2)=X, F(Y)=X\otimes X$.

{\bf Proof.} Following \cite{M2} one can describe the category $Z(\cC)$ in
the following way: the object $A=\be\boxtimes \be \oplus X\boxtimes X$ of
the category $\cC \boxtimes \cC^{op}$ has a natural structure of Frobenius 
algebra and the category $Z(\cC)$ is equivalent to the category of 
$A-$bimodules. Note that any (say left) $A-$module is free, that is of the 
form $A\otimes M$ where $M$ is some object of $\cC\boxtimes \be$ (and
$A-$module structure is the obvious one); this is a consequence of the general
fact that the category of $A-$modules considered as the module category over
$\cC \boxtimes \cC^{op}$ is module equivalent to the category $\cC$ with
$X\boxtimes Y\in \cC \boxtimes \cC^{op}$ acting via the functor $X\otimes ?
\otimes Y$, see \cite{O1}. So there are just two simple
$A-$modules --- $A$ itself and $A\otimes (X\boxtimes \be)=\be \boxtimes X+
X\boxtimes \be +nX\boxtimes X$.

We will denote by $\Hom(?,?)$ the Hom-spaces in the category $\cC \boxtimes
\cC^{op}$ and by $\Hom_{A-A}(?,?)$ the Hom-spaces in the category of 
$A-$bimodules. For any simple object $M\in \cC \boxtimes \cC^{op}$ one
considers the ``free'' bimodule $A\otimes M\otimes A$. Note that for any
$A-$bimodule $B$ one has $\Hom_{A-A}(A\otimes M\otimes A, B)=
\Hom(M, B)$. In particular any $A-$bimodule is a direct summand of some 
free bimodule. Now taking $M=\be \boxtimes \be$ we get that
$\Hom_{A-A}(A\otimes A, A\otimes A)=\Hom(\be \boxtimes \be,
A\otimes A)$ is two dimensional, so $A\otimes A$ is a direct sum of two
nonisomorphic bimodules. One of them is $A$ itself, so $A\otimes A=A\oplus Y$
where simple bimodule $Y$ as an object of $\cC$ has the following decomposition
$Y=\be \boxtimes \be \oplus n\be \boxtimes X\oplus nX\boxtimes \be \oplus
(n^2+1)X\boxtimes X$. This implies that for $M=\be \boxtimes X$ or $X\boxtimes 
\be$ one has $\Hom_{A-A}(A\otimes M\otimes A, Y)$ is $n-$dimensional
and the calculation similar to the one above shows that $A\otimes M\otimes A=
nY\oplus X_1\oplus X_2$ where $X_1, X_2$ are simple nonisomorphic bimodules
and $X_1\oplus X_2=2\be \boxtimes X\oplus 2X\boxtimes \be \oplus 2nX\boxtimes 
X$. Since $X_1$ and $X_2$ are in particular $A-$modules one gets from the 
description of $A-$modules above that $X_1=X_2=\be \boxtimes X\oplus 
X\boxtimes \be \oplus X\boxtimes X$ as objects of $\cC \boxtimes \cC^{op}$.
Finally one calculates easily that for $M=X\boxtimes X$ one has
$A\otimes M\otimes A=A\oplus nX_1\oplus nX_2\oplus (n^2+1)Y$ and so all
$A-$bimodules are classified.

The forgetful functor $F$ has the following description on the objects of
the category of $A-$bimodules: any $A-$bimodule $B$ is in particular left 
$A-$module, so is of the form $A\otimes (M\boxtimes \be)$ where $M\in \cC$;
then $F(B)=M$. This finishes the proof of the Lemma. \sq
 
The following Lemma calculates the fusion rules of $Z(\cC)$ and shows that
the based ring $K(Z(\cC))=K_n\boxtimes K_n$.

{\bf Lemma 2.2.} We have
$$X_1\otimes X_1=\be \oplus nX_1;\; X_2\otimes X_2=\be \oplus nX_2;\;
X_1\otimes X_2=X_2\otimes X_1=Y.$$

{\bf Proof.} It is known that the category $Z(\cC)$ is rigid, see \cite{M2}.
It is clear that $\be^*=\be$ and $Y^*=Y$. We claim that $X_1^*=X_1$ and
$X_2^*=X_2$. Indeed otherwise $\Hom(\be, X_1\otimes X_1)=0$ and hence 
$X_1\otimes X_1=Y$ (since $\Hom(\be,F(X_1\otimes X_1))\ne 0$. Similarly
$X_2\otimes X_2=Y$. We have also that $X_1\otimes X_2$ is a direct sum of 
$\be$ and $n$ summands each of which is isomorphic either to $X_1$ or to 
$X_2$. We can assume that $\Hom(X_1\otimes X_2,X_1)\ne 0$ (otherwise take
$X_1$ instead of $X_2$). But in this case we have a contradiction:
$$0=\Hom(X_1,Y)=\Hom(X_1,X_1\otimes X_1)=$$
$$=\Hom(X_1\otimes X_1^*,X_1)=\Hom(X_1\otimes X_2,X_1)\ne 0.$$
Thus $X_1^*=X_1$, $X_2^*=X_2$ and hence $X_1\otimes X_2=X_2\otimes X_1=Y$. 
Note that
$$\Hom(X_1\otimes X_1,X_2)=\Hom(X_1,X_1^*\otimes X_2)=\Hom(X_1,Y)=0$$
so $X_1\otimes X_1=\be \oplus nX_1$. Similarly $X_2\otimes X_2=\be \oplus 
nX_2$. The Lemma is proved. \sq

Thus the subcategory $\langle \be, X_1\rangle$ of $Z(\cC)$ consisting of 
direct sums of $\be$ and $X_1$ is a monoidal subcategory. The forgetful
functor $F$ restricted to this subcategory is an equivalence of categories
and thus we proved

{\bf Corollary 2.1.} The category $\cC$ admits a structure of braided category.
\sq

\subsection{The category $\cC$ is modular} Let us fix a structure of a braided
category on $\cC$. For $M,N\in \cC$ let $\beta_{M,N}:M\otimes N\to N\otimes M$
denote the braiding morphism. The morphisms $\beta_{M,N}$ are completely
determined by 4 morphisms $\beta_{M,N}$ where $M$ and $N$ are simple objects
of $\cC$ (since the braiding is functorial). It follows from the axioms that
$\beta_{\be,\be}=\Id$, $\beta_{\be,X}=\Id$, $\beta_{X,\be}=\Id$ so the only
nontrivial morphism is $\beta_{X,X}$. The morphism $\beta_{X,X}$ induces a
linear automorphisms of the one dimensional space $\Hom(\be,X\otimes X)$ and of
$n-$dimensional space $\Hom(X,X\otimes X)$; so the first is just some number
$\mu \in k^*$ and the second is some linear operator $\Lambda$.

{\bf Lemma 2.3.} One has $\Lambda^2=\mu \Id$.

{\bf Proof.} The vector space $\Hom(\be,X\otimes X\otimes X)$ carries the
action of two linear operators $\Lambda_1:=\beta_{X,X}\otimes \Id$ and
$\Lambda_2:=\Id \otimes \beta_{X,X}$ and it is enough to prove that 
$\Lambda_1^2=\mu \Id$. Note that by the hexagon axiom 
$(\Id \otimes \beta_{X,X})\circ (\beta_{X,X}\otimes \Id)=\beta_{X,X\otimes X}$
and hence $\Lambda_2\Lambda_1=\mu \Id$. On the other hand the braid relation
(= Yang-Baxter relation) says $\Lambda_1\Lambda_2\Lambda_1=\Lambda_2\Lambda_1
\Lambda_2$ whence $\Lambda_1=\Lambda_2$ and the Lemma is proved. \sq

{\bf Corollary 2.2.} The category $\cC$ admits a structure of ribbon 
(= pivotal and braided) category.

{\bf Proof.} It is enough to define the twists by $\theta_\be =1$ and 
$\theta_X=\mu^{-1}$ (see \cite{BK} for notations). \sq

{\bf Corollary 2.3.} Assume that $n\ne 0$. Then the category $\cC$ is modular.

{\bf Proof.} It is easy to see that if $\cC$ is not modular then $\mu =1$ and
$\cC$ is symmetric. But in a symmetric fusion category the dimensions of all
objects are integers, see \cite{AEG} Theorem 7.2. The Lemma is proved since
the quadratic equation $x^2=1+nx$ has an integer root only for $n=0$. \sq 

\subsection{The category $\cC$ does not exist} The main result of this note
is the following

{\bf Theorem 2.1.} Assume that $n\ge 2$. There is no fusion category $\cC$
such that $K(\cC)=K_n$.

{\bf Proof.} We already proved that category $\cC$ is modular if it exists. 
Let $d$ be a dimension of $X$ and let $\theta =\theta_X$. Thus $d$ is a root 
of equation $d^2=1+nd$ and $\theta$ is a root of unity by Vafa's Theorem,
see \cite{BK} Theorem 3.1.19. Consider the Gaussian sums $p_+=1+\theta d^2$ and
$p_-=1+\theta^{-1}d^2$. Since $\cC$ is a modular category one has
$p_+p_-=1+d^2$ (see \cite{BK} 3.1.15, 3.1.22) or, equivalently, $\theta +
\theta^{-1}=nd$. We can assume that $d>n$ (applying otherwise the Galois 
automorphism to the equation). Then $|nd|>n^2$ and 
$|\theta +\theta^{-1}|\le 2$ and we
have a contradiction for $n\ge 2$. \sq

\subsection{Categories $\cC$ with $K(\cC)=K_0$} Let $G$ be a finite group
and let $K(G)$ be the based ring with basis $X_g, g\in G$ and fusion rules
$X_g\otimes X_h=X_{gh}$. It is shown in \cite{MS} Appendix E that the
monoidal categories $\cC$ with $K(\cC)=K(G)$ are classified by $H^3(G,k^*)$.
Any category of this kind is automatically rigid.
In particular $K_0=K(\BZ/2\BZ)$ and categorifications of $K_0$ are
classified by $H^3(\BZ/2\BZ ,k^*)=\BZ/2\BZ$. So there are two such categories.
First one is the category of representations of $\BZ/2\BZ$ and second differs
from the first one by the sign of associativity morphism $(X\otimes X)\otimes X
\to X\otimes (X\otimes X)$. This second category can be explicitly realized
as the fusion category of integrable representations of $\widehat{sl}_2$ at
level 1.

\subsection{Categories $\cC$ with $K(\cC)=K_1$} Let $\cC$ be a monoidal
category with $K(\cC)=K_1$. Choose basis vectors in one 
dimensional vector spaces $\Hom(\be,X\otimes X)$ and $\Hom(X,X\otimes X)$. 
The only nontrivial associativity morphism is $(X\otimes X)\otimes X\to
X\otimes (X\otimes X)$, this is equivalent to giving isomorphisms of
vector spaces $\Hom(X,X\otimes X)\otimes \Hom(\be,X\otimes X)\to
\Hom(X,X\otimes X)\otimes \Hom(\be,X\otimes X)$ (so this is just a number
$\lambda \in k^*$) and $\Hom(\be,X\otimes X)\otimes \Hom(X,X\otimes \be)\oplus
\Hom(X,X\otimes X)\otimes \Hom(X,X\otimes X)\to \Hom(\be,X\otimes X)\otimes 
\Hom(X,\be \otimes X)\oplus \Hom(X,X\otimes X)\otimes \Hom(X,X\otimes X)$
(so this can be represented by invertible $2\times 2$ matrix). After relatively
easy calculation (along the lines of \cite{TY}) one finds that $\lambda =1$
and the second isomorphism can be represented by the matrix $\left(
\begin{array}{cc}a&1\\a&-a\end{array}\right)$ where $a$ is a root of equation
$a^2+a=1$. So we have two solutions; both categories are rigid; they differ
by the action of Galois group. In one of these categories $\dim(X)=\frac{1+
\sqrt{5}}{2}$ and in the second $\dim(X)=\frac{1-\sqrt{5}}{2}$. The first
category can be explicitly realized as a subcategory of ``integer spin''
representations of the fusion category of integrable $\widehat{sl}_2-$modules
at level 3 and the second category is the minimal model $\cM (2,5)$ for
the Virasoro algebra (with central charge $c=-\frac{22}{5}$). Both categories
can be also realized using the quantum group $U_q(sl_2)$ for $q=\sqrt[10]{1}$
(there are 4 primitive tenth roots of 1, but $q$ and $q^{-1}$ give rise to
the same category), see \cite{BK}.

\section{Appendix}
\centerline{Pavel Etingof}
\bigskip

In this appendix we give an upper bound for the dimension of a semisimple
Hopf algebra over $k$ with $n$ irreducible representations. 
For group algebras, this estimate was obtained 100 years ago by
E. Landau, \cite{L}. (We are grateful to R.Stanley for giving us
this reference).

For a positive integer $n$, let $P(n)$ be the set of 
positive integer solutions of the equation
$\frac{1}{x_1}+...+\frac{1}{x_n}=1$, and let 
$d(n)={\rm max}_{(x_1,...,x_n)\in P(n)}{\rm max}_i x_i$. 

{\bf Theorem 3.1.} Let $H$ be a semisimple Hopf algebra over $k$,
which has $n$ irreducible representations. Then
${\rm dim} H\le d(n)$. 

{\bf Remark 1.} The set $P(n)$ is finite, since for any 
positive rational number $r$ the number of positive integer
solutions of the equation $\frac{1}{x_1}+...+\frac{1}{x_n}=r$
is finite (\cite{L}). Indeed, if $x_1$ is the smallest
coordinate, then $x_1\le n/r$, so there are finitely many
possibilities for $x_1$, and if $x_1$ is fixed then 
$(x_2,...,x_n)$ vary over the set of positive integer solutions of 
the equation $\frac{1}{x_2}+...+\frac{1}{x_n}=r-\frac{1}{x_1}$,  
so the statement follows by induction. (In fact, this inductive
procedure allows one to obtain an explicit estimate for $d(n)$, see
\cite{L}). 

{\bf Remark 2.} One has $d(n+1)\ge 2d(n)$ since 
if $\sum_{i=1}^n\frac{1}{x_i}=1$ then 
$\sum_{i=1}^n\frac{1}{2x_i}+\frac{1}{2}=1$. 

\begin{proof}
Our proof is a generalization of the classical proof of Landau in
the group case, 
which employs the fact 
that for a group with $n$ irreducible representations, one has 
$\sum_{i=1}^n\frac{1}{|C_i|}=1$, where $C_i$ are centralizers of
conjugacy classes (the class equation). Namely, we 
use the Hopf algebraic class equation, due to G.Kac
\cite{Ka} and Y.Zhu \cite{Zhu}.

Let $C(H)\subset H^*$ be the character ring of $H$ (it is spanned
by characters of irreducible representation). This is a
semisimple algebra of dimension $n$, so $C(H)$ can be
identified with $\oplus {\rm
  Mat}_{r_i}$, where $\sum r_i^2=n$. Let us choose such an
identification,
and let $E_{jj}^{(i)}\in C(H)$ be the corresponding matrix
units. 

{\bf Theorem} (Kac-Zhu), \cite{Zhu} The ratio $m_i={\rm dim(H)}/{\rm
  Tr}|_{H^*}(E_{jj}^{(i)}\cdot)$ (where for $a\in H^*$ $a\cdot$ denotes the
operator of multiplication by $a$) is a positive integer.

Now observe that $\sum_i\sum_{j=1}^{r_i}E_{jj}^{i}=1$, so 
$\sum_i\sum_{j=1}^{r_i}\frac{r_i}{m_i}=1$. Writing 
$r_i/m_i$ as a sum of $r_i$ copies of $1/m_i$, we get a solution
of the equation $\sum_{j=1}^r\frac{1}{x_i}=1$, where $r=\sum_i
r_i\le n$. Thus, $m_i\le d(r)\le d(n)$. On the other hand, 
consider the 1-dimensional matrix block of $C(H)$ spanned by the
integral of $H^*$ (i.e. by the character of the regular
representation of $H$). For this block the number $m_i$ is obviously
equal to ${\rm dim}(H)$, since the integral is a projector to a
1-dimensional subspace. So ${\rm dim}H\le d(n)$, as desired. 
\end{proof}

\end{document}